\documentclass[a4paper,11pt]{article}

\usepackage{amsfonts, amssymb, amscd, amsmath, enumerate, verbatim, amsthm}
\usepackage{array}
\usepackage{changepage}
\usepackage[margin=1.4in]{geometry}
\usepackage[T1]{fontenc}
\usepackage{algorithm,algorithmic}
\usepackage{multirow}
\usepackage{tikz}

\usepackage{hyperref}

\begin{document}
	
	\title{A Short Survey of Averaging Techniques in Stochastic Gradient Methods}
	\author{K. Lakshmanan \\ Department of Computer Science and Engineering \\ Indian Insititute of Technology (BHU), Varanasi 221005.\\ Email: lakshmanank.cse@iitbhu.ac.in}
	\date{}
	
	\maketitle

	\begin{abstract}
		
		Stochastic gradient methods are among the most widely used algorithms
		for large-scale optimization and machine learning. A key technique for
		improving the statistical efficiency and stability of these methods is
		the use of averaging schemes applied to the sequence of iterates
		generated during optimization. Starting from the classical work on
		stochastic approximation, averaging techniques such as
		Polyak--Ruppert averaging have been shown to achieve optimal
		asymptotic variance and improved convergence behavior.
		
		In recent years, averaging methods have gained renewed attention in
		machine learning applications, particularly in the training of deep
		neural networks and large-scale learning systems. Techniques such as
		tail averaging, exponential moving averages, and stochastic weight
		averaging have demonstrated strong empirical performance and improved
		generalization properties.
		
		This paper provides a survey of averaging techniques in stochastic
		gradient optimization. We review the theoretical foundations of
		averaged stochastic approximation, discuss modern developments in
		stochastic gradient methods, and examine applications of averaging
		in machine learning. In addition, we summarize recent results on the
		finite-sample behavior of averaging schemes and highlight several
		open problems and directions for future research.	
	\end{abstract}

\section{Introduction}

Stochastic optimization methods have become a fundamental tool in modern
machine learning and large-scale data analysis. Many learning problems
can be formulated as minimizing an expected loss function of the form
\[
\min_{x \in \mathbb{R}^d} \; f(x) = \mathbb{E}_{\xi}[F(x,\xi)],
\]
where $\xi$ denotes a random data sample. In such settings, stochastic
approximation methods estimate gradients using randomly sampled data
points and update parameters iteratively. The classical stochastic
approximation framework originates from the seminal work of
Robbins and Monro~\cite{RobbinsMonro1951}, with further developments
by Kiefer and Wolfowitz~\cite{KieferWolfowitz1952} and later analyses
of asymptotic properties such as those of Fabian~\cite{Fabian1968}.
Over the decades, stochastic gradient descent (SGD) has emerged as
one of the most widely used algorithms for solving large-scale
optimization problems, particularly in machine learning applications
where computing full gradients is computationally prohibitive.

A key development in stochastic approximation is the idea of
\emph{averaging the iterates} generated by SGD. Early work by
Polyak~\cite{Polyak1990} and the influential analysis of
Polyak and Juditsky~\cite{PolyakJuditsky1992} showed that averaging
stochastic gradient iterates can significantly improve statistical
efficiency and achieve optimal asymptotic variance under suitable
conditions. Related contributions by Ruppert~\cite{Ruppert1988}
and later developments in stochastic optimization theory further
established averaged stochastic gradient methods as a central
technique in both statistics and machine learning. In recent years,
averaging techniques have also gained renewed interest in deep
learning, where methods such as stochastic weight averaging and
model averaging have been shown to improve generalization
performance and stability of training~\cite{Izmailov2018}.

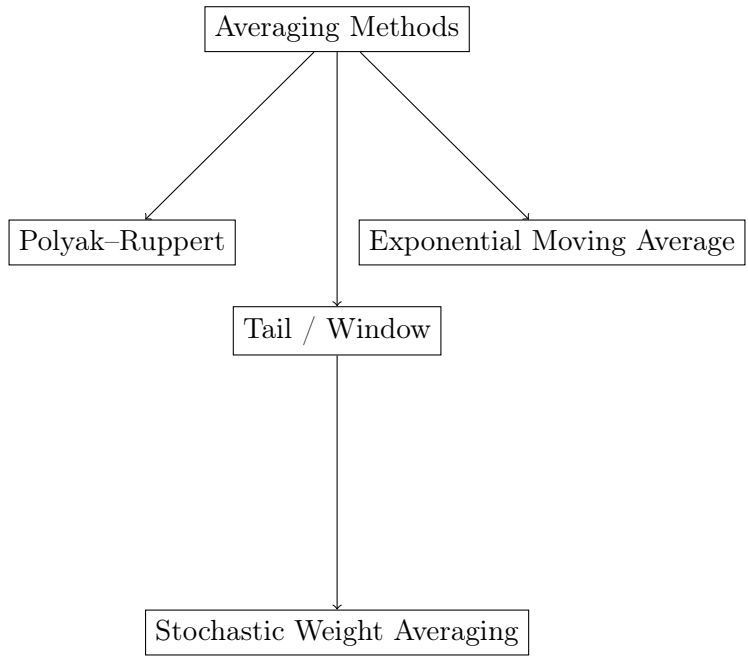
\begin{figure}[t]
	\centering
	\begin{tikzpicture}[node distance=4cm]
		
		\node (avg) [draw, rectangle] {Averaging Methods};
		
		\node (polyak) [draw, rectangle, below left of=avg] {Polyak--Ruppert};
		\node (tail) [draw, rectangle, below of=avg] {Tail / Window};
		\node (ema) [draw, rectangle, below right of=avg] {Exponential Moving Average};
		
		\node (swa) [draw, rectangle, below of=tail] {Stochastic Weight Averaging};
		
		\draw[->] (avg) -- (polyak);
		\draw[->] (avg) -- (tail);
		\draw[->] (avg) -- (ema);
		\draw[->] (tail) -- (swa);
		
	\end{tikzpicture}
	\caption{Main categories of averaging techniques in stochastic optimization.}
\end{figure}

Despite the widespread use of averaging in stochastic optimization,
the literature on the subject is scattered across several research
communities, including statistics, optimization, and machine learning.
Early work focused primarily on asymptotic efficiency and statistical
properties of averaged stochastic approximation algorithms. More
recent studies have investigated finite-sample convergence rates,
robustness properties, and practical implementations of averaging
schemes in large-scale learning systems. In parallel, deep learning
research has introduced new forms of averaging, such as exponential
moving averages and weight averaging techniques, which have been
empirically shown to improve generalization and training stability.
These developments highlight the need for a unified overview of
averaging techniques and their theoretical foundations.

The goal of this survey is to provide a comprehensive overview of
averaging methods in stochastic gradient algorithms. We review the
historical development of averaging in stochastic approximation,
summarize theoretical results on convergence and statistical
efficiency, and discuss modern developments in machine learning
applications. In particular, we examine how averaging interacts with
stochastic noise, step-size schedules, and problem structure, and
how these factors influence the performance of stochastic gradient
methods in practice.

\begin{table}[t]
	\setlength{\tabcolsep}{0pt}
	\centering
	\caption{Overview of averaging techniques in stochastic optimization}
	\vspace{0.2cm}
	\centering{\parindent=-0.5em\scriptsize
	\begin{tabular}{|l|l|l|l|}
		\hline
		Method & Averaging Rule & Key Idea & Reference \\
		\hline
		Polyak--Ruppert Averaging & $\bar{x}_k=\frac{1}{k}\sum_{i=1}^{k}x_i$ & Uniform average of all iterates & \cite{PolyakJuditsky1992,Ruppert1988} \\
		\hline
		Tail Averaging & $\bar{x}_k=\frac{1}{m}\sum_{i=k-m+1}^{k}x_i$ & Average recent iterates & \cite{ShamirZhang2013} \\
		\hline
		Window Averaging & Sliding window average & Discard early transient iterates & Various \\
		\hline
		Weighted Averaging & $\bar{x}_k=\sum_{i=1}^{k}w_i x_i$ & Non-uniform weights & Various \\
		\hline
		Exponential Moving Average & $\bar{x}_k=\beta\bar{x}_{k-1}+(1-\beta)x_k$ & Emphasize recent iterates & \cite{KingmaBa2015} \\
		\hline
		Stochastic Weight Averaging & Average selected SGD iterates & Locate wide optima & \cite{Izmailov2018} \\
		\hline
		Snapshot Ensembles & Average multiple models & Ensemble learning & \cite{Huang2017} \\
		\hline
	\end{tabular}}
\end{table}

The remainder of this paper is organized as follows.
Section~\ref{sec:background} reviews the classical framework of
stochastic approximation and stochastic gradient methods.
Section~\ref{sec:polyak} discusses the Polyak--Ruppert averaging
scheme and its theoretical properties. Section~\ref{sec:tail}
covers alternative averaging strategies such as tail averaging
and window averaging. Section~\ref{sec:ema} reviews exponential
moving averages and related techniques widely used in modern
machine learning. Section~\ref{sec:ml} discusses applications of
averaging methods in deep learning and large-scale training.
Section~\ref{sec:finite} summarizes recent results on finite-sample
behavior and theoretical guarantees for averaged stochastic
optimization. Then guidelines for practioners is given in Section~\ref{sec:guidelines}. Finally, Section~\ref{sec:open} outlines open
problems and future research directions.

\section{Background: Stochastic Approximation and Gradient Methods}
\label{sec:background}

\begin{figure}[h]
	\centering
	\begin{tikzpicture}[scale=1]
		
		\draw[->, thick] (0,0) -- (11,0);
		
		\node at (0,0) [below] {1950};
		\node at (3,0) [below] {1990};
		\node at (6,0) [below] {2010};
		\node at (9,0) [below] {2020};
		
		\draw (1,0.1) -- (1,-0.1);
		\node[align=center] at (1,1) {Robbins--Monro\\Stochastic Approximation\\(1951)};
		
		\draw (3.5,0.1) -- (3.5,-0.1);
		\node[align=center] at (3.5,1) {Polyak--Ruppert\\Averaging\\(1992)};
		
		\draw (6.5,0.1) -- (6.5,-0.1);
		\node[align=center] at (6.5,1) {Modern SGD\\Theory\\(2010s)};
		
		\draw (9.5,0.1) -- (9.5,-0.1);
		\node[align=center] at (9.5,1) {Stochastic Weight\\Averaging\\(2018)};
		
	\end{tikzpicture}
	
	\caption{Historical development of averaging techniques in stochastic optimization.}
\end{figure}
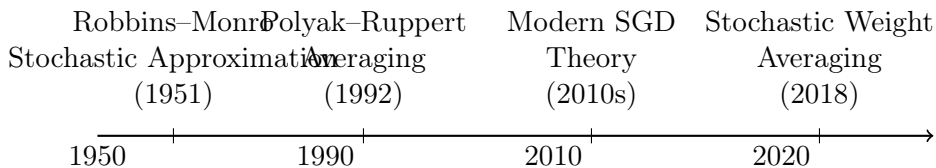

Stochastic approximation methods provide a general framework for solving
optimization and root-finding problems when only noisy observations of
the objective or gradient are available. These methods arise naturally
in machine learning, signal processing, and statistics, where the
objective function often involves expectations over random data samples.

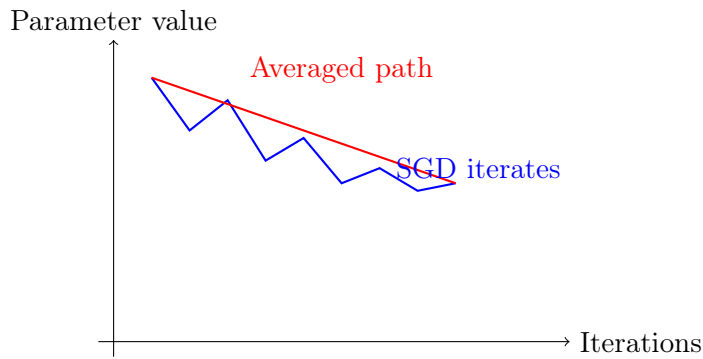
\begin{figure}[h]
	\centering
	\begin{tikzpicture}[scale=1]
		
		\draw[->] (-0.2,0) -- (6,0) node[right] {Iterations};
		\draw[->] (0,-0.2) -- (0,4) node[above] {Parameter value};
		
		\draw[blue, thick] 
		(0.5,3.5) -- (1,2.8) -- (1.5,3.2) -- (2,2.4) -- (2.5,2.7) -- (3,2.1)
		-- (3.5,2.3) -- (4,2.0) -- (4.5,2.1);
		
		\draw[red, thick] (0.5,3.5) -- (4.5,2.1);
		
		\node[blue] at (4.8,2.3) {SGD iterates};
		\node[red] at (3,3.6) {Averaged path};
		
	\end{tikzpicture}
	\caption{Stochastic gradient iterates exhibit noisy behavior,
		while averaged iterates provide a smoother trajectory.}
\end{figure}

\subsection{Stochastic Approximation}

Consider the problem of finding the root of a function
\[
h(x) = \mathbb{E}_{\xi}[H(x,\xi)],
\]
where $\xi$ is a random variable and only noisy observations
$H(x,\xi)$ are available. The classical stochastic approximation
algorithm introduced by Robbins and Monro~\cite{RobbinsMonro1951}
generates a sequence $\{x_k\}$ according to the recursion
\[
x_{k+1} = x_k - \eta_k H(x_k,\xi_k),
\]
where $\eta_k > 0$ is a step-size sequence and $\{\xi_k\}$ are
independent samples. Under suitable conditions on the step sizes
and regularity assumptions on the function $h$, the sequence
$\{x_k\}$ converges to a solution of the equation $h(x)=0$.

The stochastic approximation framework has been extensively studied
in the literature. Early developments include the work of
Kiefer and Wolfowitz~\cite{KieferWolfowitz1952}, who considered
gradient-free stochastic approximation methods, and the asymptotic
analysis of Fabian~\cite{Fabian1968}. Modern treatments of stochastic
approximation can be found in the monographs of
Kushner and Yin~\cite{KushnerYin2003} and
Borkar~\cite{Borkar2008}.

\subsection{Stochastic Gradient Descent}

In many machine learning applications the objective function can be
written as an expectation
\[
f(x) = \mathbb{E}_{\xi}[F(x,\xi)],
\]
where $\xi$ represents a random data point and $F(x,\xi)$ denotes the
loss associated with parameter $x$ and sample $\xi$.
The gradient of $f$ is then given by
\[
\nabla f(x) = \mathbb{E}_{\xi}[\nabla F(x,\xi)] .
\]

Computing the exact gradient is often computationally expensive when
the dataset is large \cite{la2025gradient}. Instead, stochastic gradient descent (SGD)
uses a noisy estimate of the gradient based on a randomly sampled
data point. The SGD iteration takes the form
\[
x_{k+1} = x_k - \eta_k g_k,
\]
where $g_k = \nabla F(x_k,\xi_k)$ is an unbiased estimate of the
true gradient $\nabla f(x_k)$ and $\eta_k$ is the learning rate
or step size.

SGD has become one of the most widely used optimization methods in
machine learning due to its simplicity and scalability.
A comprehensive overview of optimization methods used in
large-scale machine learning can be found in
Bottou et al.~\cite{Bottou2018}.

\subsection{Convergence Properties}

The convergence behavior of stochastic gradient methods depends
critically on the choice of step-size sequence $\{\eta_k\}$.
For classical stochastic approximation \cite{Borkar2008}, it is typically assumed
that the step sizes satisfy
\[
\sum_{k=1}^{\infty} \eta_k = \infty,
\qquad
\sum_{k=1}^{\infty} \eta_k^2 < \infty .
\]

Under suitable smoothness and convexity assumptions on the objective
function, these conditions ensure convergence of the iterates to a
stationary point. However, the stochastic noise inherent in gradient
estimates often leads to slow convergence and large variance in the
iterates. As a result, techniques for reducing variance and improving
the statistical efficiency of SGD have become an important topic of
research in stochastic optimization.

One of the most influential approaches for reducing the variance of
stochastic gradient iterates is the use of \emph{iterate averaging}.
Instead of returning the final iterate $x_k$, averaged stochastic
approximation methods compute an average of past iterates,
\[
\bar{x}_k = \frac{1}{k} \sum_{i=1}^{k} x_i .
\]
As shown in the work of Polyak and Juditsky~\cite{PolyakJuditsky1992},
such averaging can significantly improve the statistical efficiency
of stochastic approximation algorithms and lead to optimal asymptotic
performance.

In the following sections we review the development of averaging
methods for stochastic gradient algorithms, starting with the
classical Polyak--Ruppert averaging scheme and its theoretical
properties.

\section{Polyak--Ruppert Averaging}
\label{sec:polyak}

One of the most influential developments in stochastic approximation
is the idea of averaging the iterates produced by stochastic gradient
methods. Instead of returning the final iterate of a stochastic
approximation algorithm, Polyak and Juditsky~\cite{PolyakJuditsky1992}
proposed returning the average of all previous iterates. This simple
modification significantly improves the statistical efficiency of the
algorithm and has become a central technique in stochastic optimization.

\subsection{The Averaged Stochastic Gradient Algorithm}

Consider the stochastic gradient descent iteration
\[
x_{k+1} = x_k - \eta_k g_k,
\]
where $g_k$ is an unbiased estimator of the gradient
$\nabla f(x_k)$ and $\eta_k$ is a step size.
The averaged iterate is defined as
\[
\bar{x}_k = \frac{1}{k} \sum_{i=1}^{k} x_i .
\]

The key idea of Polyak--Ruppert averaging \cite{PolyakRuppertReview} is that the sequence
$\{\bar{x}_k\}$ often converges faster and exhibits smaller variance
than the sequence of individual iterates $\{x_k\}$. In particular,
the averaged estimator effectively smooths out the noise introduced
by stochastic gradients.

\subsection{Asymptotic Optimality}

A fundamental result established by Polyak and Juditsky
\cite{PolyakJuditsky1992} shows that averaged stochastic approximation
achieves optimal asymptotic variance under suitable regularity
conditions. Suppose the objective function $f$ has a unique minimizer
$x^*$ and the stochastic gradient satisfies
\[
\mathbb{E}[g_k \mid x_k] = \nabla f(x_k).
\]
Under appropriate assumptions on the step sizes and smoothness of the
objective function, the averaged estimator satisfies
\[
\sqrt{k}(\bar{x}_k - x^*) \;\xrightarrow{d}\; 
\mathcal{N}(0, \Sigma),
\]
where $\Sigma$ is the optimal covariance matrix associated with the
underlying stochastic approximation process.

This result shows that averaging recovers the statistical efficiency
of the optimal estimator, even when the step sizes used in the
stochastic gradient method are relatively large. As a result,
Polyak--Ruppert averaging allows stochastic gradient methods to
achieve the same asymptotic efficiency as classical statistical
estimators.

\subsection{Historical Development}

The idea of iterate averaging was originally proposed by
Polyak~\cite{Polyak1990} and further analyzed in the work of
Polyak and Juditsky~\cite{PolyakJuditsky1992}. Independently,
Ruppert~\cite{Ruppert1988} studied similar averaging schemes in
the context of stochastic approximation and statistical estimation.
These works demonstrated that averaging could dramatically improve
the performance of stochastic approximation algorithms.

Subsequent research extended these results to broader classes of
optimization problems. Nemirovski et al.~\cite{Nemirovski2009}
analyzed averaged stochastic gradient methods for stochastic
programming problems and showed improved convergence guarantees
in convex optimization settings. Later work by Bach and
Moulines~\cite{BachMoulines2013} provided non-asymptotic convergence
rates for averaged stochastic gradient methods, showing that
averaging achieves optimal convergence rates in certain convex
optimization problems.

\subsection{Interpretation as Variance Reduction}

One way to interpret Polyak--Ruppert averaging is as a form of
variance reduction. Individual stochastic gradient iterates
often exhibit large fluctuations due to the randomness of gradient
estimates. By averaging multiple iterates, the variance of the
estimator is reduced, leading to more stable convergence.

From a statistical perspective, the averaged iterate can be viewed
as an estimator constructed from multiple noisy observations of the
optimization process. The averaging procedure effectively filters
out high-frequency noise in the iterates, resulting in improved
estimation accuracy.

\paragraph{Limitations and Practical Considerations}

Although Polyak--Ruppert averaging has strong theoretical guarantees,
its practical performance depends on several factors. For example,
averaging all iterates may include early iterations that are far
from the optimum, potentially slowing convergence in practice.
This observation has motivated alternative averaging schemes,
such as tail averaging and window averaging, which average only
a subset of the most recent iterates.

These alternative strategies have been explored in both theoretical
and empirical studies of stochastic optimization algorithms.
In the next section we review several such variants and discuss
their advantages and limitations.

\section{Tail Averaging and Window Averaging}
\label{sec:tail}

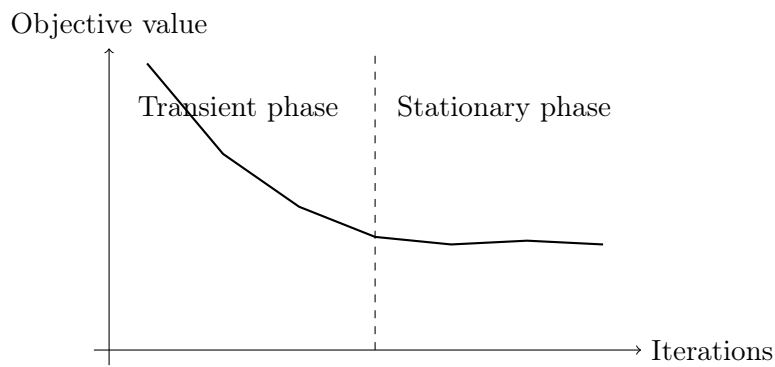
\begin{figure}[h]
	\centering
	\begin{tikzpicture}[scale=1]
		
		\draw[->] (-0.2,0) -- (7,0) node[right] {Iterations};
		\draw[->] (0,-0.2) -- (0,4) node[above] {Objective value};
		
		\draw[thick] 
		(0.5,3.8) -- (1.5,2.6) -- (2.5,1.9) -- (3.5,1.5)
		-- (4.5,1.4) -- (5.5,1.45) -- (6.5,1.4);
		
		\draw[dashed] (3.5,0) -- (3.5,4);
		
		\node at (1.7,3.2) {Transient phase};
		\node at (5.2,3.2) {Stationary phase};
		
	\end{tikzpicture}
	\caption{Stochastic optimization typically exhibits a transient phase
		followed by a stationary phase where iterates fluctuate near the optimum.}
\end{figure}

Although Polyak--Ruppert averaging provides strong theoretical
guarantees, averaging all iterates generated by stochastic gradient
descent may not always be optimal in practice. In particular,
early iterates of stochastic optimization algorithms are often
far from the optimal solution and may introduce unnecessary bias
into the averaged estimator. This observation has motivated several
variants of averaging schemes that focus on averaging only the
most recent iterates of the optimization process.

\subsection{Tail Averaging}

A commonly used alternative to full averaging is \emph{tail
	averaging}. Instead of averaging all iterates $\{x_1,\ldots,x_k\}$,
tail averaging computes the average of the last $m$ iterates,
\[
\bar{x}_k = \frac{1}{m} \sum_{i=k-m+1}^{k} x_i .
\]

The motivation for tail averaging is that stochastic gradient
iterates typically go through an initial \emph{transient phase}
during which the algorithm moves rapidly toward the optimum.
During this phase the iterates may be highly biased and averaging
them can degrade performance. Once the algorithm enters a
\emph{stationary regime}, the iterates fluctuate around the
optimal solution due to stochastic noise. Averaging in this
stationary phase can reduce variance and improve stability.

Theoretical analyses of tail averaging have shown that it can
retain many of the favorable properties of Polyak--Ruppert
averaging while improving finite-sample performance. In
particular, Shamir and Zhang~\cite{ShamirZhang2013} studied
optimal averaging schemes for stochastic gradient methods
in nonsmooth optimization and demonstrated that averaging
strategies focusing on later iterates can achieve improved
convergence behavior.

\subsection{Window Averaging}

A related idea is \emph{window averaging} \cite{fu2014stochastic}, in which the average
is computed over a fixed-size sliding window of recent iterates.
Given a window size $m$, the averaged iterate at iteration $k$
is defined as
\[
\bar{x}_k = \frac{1}{m} \sum_{i=k-m+1}^{k} x_i .
\]

Window averaging allows the algorithm to maintain a running
average that adapts as the optimization proceeds. Unlike full
averaging, which assigns equal weight to all past iterates,
window averaging focuses on the most recent iterates that are
likely to be closer to the optimum.

From a practical perspective, window averaging also has
computational advantages, since it requires storing only a
limited number of past iterates rather than the entire
optimization trajectory.

\subsection{Weighted Averaging Schemes}

More general averaging schemes assign weights to past iterates \cite{dippon1997weighted}
according to a specified weighting function. A general weighted
average takes the form
\[
\bar{x}_k = \sum_{i=1}^{k} w_i x_i ,
\]
where the weights satisfy $w_i \ge 0$ and
$\sum_{i=1}^{k} w_i = 1$.

Several weighting strategies have been proposed in the
literature. For example, exponentially decaying weights
emphasize recent iterates and lead to the well-known
\emph{exponential moving average} technique widely used
in machine learning algorithms. Other schemes assign
weights that increase with iteration number, giving
greater importance to later iterates that are closer to
the optimum.

These weighted averaging methods provide additional
flexibility in controlling the bias–variance trade-off
in stochastic optimization. In particular, assigning
greater weight to later iterates can help reduce bias
introduced during the transient phase of the algorithm.

\paragraph{Practical Considerations}

In practical machine learning applications, the choice
of averaging strategy often depends on the characteristics
of the optimization problem and the training procedure.
Full averaging is attractive from a theoretical perspective,
but tail or window averaging is often preferred in large-scale
learning systems where the early optimization trajectory
may be highly nonstationary.

Recent empirical studies in deep learning have also shown
that averaging techniques can improve generalization
performance of neural networks. In particular, methods
based on averaging model parameters across multiple
training iterations have been shown to produce models
that lie in flatter regions of the loss landscape and
therefore generalize better~\cite{Izmailov2018}.

These observations have motivated further research into
adaptive and problem-dependent averaging schemes for
stochastic optimization.

\section{Averaging in Modern Machine Learning}
\label{sec:ml}

In recent years averaging techniques have become widely used in
machine learning \cite{ShamirZhang2013,Gower2019}, particularly in the training of large-scale
models such as deep neural networks. While classical stochastic
approximation literature focused primarily on statistical efficiency
and asymptotic properties \cite{Moulines2011}, modern machine learning applications
emphasize generalization performance, training stability, and
robustness to optimization noise. As a result, several averaging
methods have been proposed and adopted in practical learning systems.

\subsection{Exponential Moving Averages}

One of the most commonly used averaging techniques in machine
learning is the \emph{exponential moving average} (EMA) of model
parameters. Given a sequence of iterates $\{x_k\}$ produced by an
optimization algorithm, the exponential moving average is defined
recursively as
\[
\bar{x}_k = \beta \bar{x}_{k-1} + (1-\beta)x_k ,
\]
where $\beta \in (0,1)$ controls the decay rate of past information.

EMA places greater weight on recent iterates while retaining
information from earlier steps of the optimization trajectory.
This method is widely used in deep learning training procedures
and is often employed to stabilize optimization and improve
generalization. Similar ideas appear in adaptive optimization
methods such as Adam~\cite{KingmaBa2015}, where exponential
averaging is used to estimate gradient moments during training.

\subsection{Stochastic Weight Averaging}

A particularly influential averaging technique in deep learning
is \emph{stochastic weight averaging} (SWA), introduced by
Izmailov et al.~\cite{Izmailov2018}. The key idea is to average
model parameters obtained at different stages of stochastic
gradient training, typically using a cyclic or constant learning
rate schedule. The averaged parameters are given by
\[
\bar{x}_K = \frac{1}{K} \sum_{k=1}^{K} x_{t_k},
\]
where $\{x_{t_k}\}$ are selected iterates from the training
trajectory.

Empirical studies show that SWA often leads to solutions that
lie in wider regions of the loss landscape. Such solutions
are associated with improved generalization performance in
deep neural networks. This observation has stimulated extensive
research into the geometry of loss surfaces \cite{Garipov2018} and the role of
averaging in locating flat minima.

\subsection{Model Averaging and Ensembles}

Another important application of averaging in machine learning
is \emph{model averaging} or ensemble learning. In this approach,
multiple models trained during different stages of optimization
are combined to produce a single prediction model. For example,
snapshot ensemble methods~\cite{Huang2017} generate several
models during a single training run by periodically varying
the learning rate and then averaging their predictions.

Model averaging can significantly improve predictive performance
by reducing variance \cite{JohnsonZhang2013} and exploiting diversity among models.
Although conceptually different from parameter averaging,
both approaches rely on similar principles of combining
multiple estimators to improve stability and accuracy.

\subsection{Distributed and Federated Learning}

Averaging also plays a central role in distributed optimization
and federated learning systems \cite{zhang2021survey}. In these settings, multiple
computing nodes or devices perform local optimization steps
and periodically communicate updates to a central server.
The server aggregates the local models using averaging
operations to produce a global model.

This paradigm, often referred to as \emph{model averaging} \cite{claeskens2008model},
forms the basis of many distributed learning algorithms.
The effectiveness of such methods depends on the interaction
between local stochastic optimization and global averaging
operations. As large-scale machine learning systems continue
to grow, understanding the role of averaging in distributed
optimization has become an active area of research.

\subsection{Generalization and Flat Minima}

Recent research suggests that averaging techniques may also
improve the generalization properties of machine learning
models. Empirical studies indicate that averaging parameters
across different stages of stochastic gradient descent tends
to produce solutions located in flatter regions of the loss
landscape. Such solutions are believed to generalize better
to unseen data.

These observations have motivated theoretical investigations
into the relationship between stochastic optimization,
loss surface geometry \cite{Garipov2018}, and averaging methods \cite{ShamirZhang2013}. Although a
complete theoretical explanation remains an open problem,
the growing body of empirical evidence suggests that
averaging plays a crucial role in the success of modern
stochastic optimization methods in machine learning.

\section{Finite-Sample Behavior of Averaging Methods}
\label{sec:finite}

\begin{table}[t]
	\centering
	\caption{Theoretical characteristics of averaging methods}
	\centering
	\vspace{0.2cm}
	{\parindent=-0.5em \scriptsize
	\begin{tabular}{|l|l|l|l|}
		\hline
		Method & Variance Reduction & Finite Sample Behavior & Typical Applications \\
		\hline
		Polyak--Ruppert Averaging & Strong & Good asymptotically & Stochastic approximation \\
		\hline
		Tail Averaging & Moderate & Improved finite sample & SGD optimization \\
		\hline
		Window Averaging & Moderate & Depends on window size & Online learning \\
		\hline
		Weighted Averaging & Adjustable & Problem dependent & Adaptive optimization \\
		\hline
		EMA & Moderate & Stable training & Deep learning \\
		\hline
		SWA & Strong empirically & Good generalization & Neural network training \\
		\hline
	\end{tabular}}
\end{table}

While classical analyses of stochastic approximation focus on
asymptotic convergence properties, modern machine learning
applications often operate in regimes where only a finite number
of optimization steps can be performed. As a result, understanding
the \emph{finite-sample behavior} of stochastic gradient methods
and their averaging variants has become an important topic of
research.

\subsection{Non-Asymptotic Convergence Rates}

Early analyses of Polyak--Ruppert averaging established
asymptotic optimality of the averaged estimator \cite{Polyak1990,PolyakRuppertReview}. However,
these results do not directly characterize the behavior
of stochastic gradient algorithms during the early stages
of optimization.

More recent work has focused on deriving non-asymptotic
convergence guarantees for stochastic gradient methods.
For example, Bach and Moulines~\cite{BachMoulines2013}
analyzed averaged stochastic gradient methods in convex
optimization settings and showed that the averaged iterates
can achieve an optimal convergence rate of order $O(1/n)$
for smooth convex problems.

Further developments in stochastic optimization theory
have provided refined analyses of stochastic gradient
methods under various structural assumptions on the
objective function, including strong convexity,
smoothness, and noise conditions. These results
provide a more detailed understanding of the behavior
of averaged stochastic gradient algorithms in finite
time.

\subsection{Bias--Variance Trade-offs}

Averaging schemes often involve a trade-off between bias
and variance in the resulting estimator. Full averaging
over all iterates can effectively reduce variance \cite{JohnsonZhang2013} but may
introduce bias due to the inclusion of early iterates
that are far from the optimal solution.

Tail averaging and window averaging schemes attempt to
reduce this bias by focusing on iterates generated after
the optimization process has approached the stationary
region near the optimum. By discarding early iterations,
these methods can achieve improved finite-sample
performance while still benefiting from variance
reduction.

Understanding the optimal balance between bias and
variance remains an active area of research. In
practice, the best averaging strategy may depend on
problem-specific factors such as noise levels,
learning rate schedules, and the geometry of the
objective function.

\subsection{Optimal Weighting Strategies}

Beyond simple averaging schemes, several studies have
investigated the use of weighted averages of stochastic
gradient iterates. In a general weighted averaging
framework, the estimator is given by
\[
\bar{x}_k = \sum_{i=1}^{k} w_i x_i,
\]
where the weights $\{w_i\}$ satisfy
$\sum_{i=1}^{k} w_i = 1$.

The choice of weights can significantly influence the
performance of the resulting estimator. In particular,
assigning larger weights to later iterates may reduce
bias arising from early transient behavior of the
optimization process. Conversely, uniform weights
provide stronger variance reduction when the iterates
are already close to the optimum.

Recent research has explored the design of weighting
schemes that adapt to the dynamics of stochastic
optimization algorithms. These approaches aim to
improve finite-sample performance while retaining the
theoretical guarantees of classical averaging methods.

\paragraph{Practical Implications}

Finite-sample analyses highlight the importance of
carefully choosing averaging strategies in practical
machine learning applications. In large-scale learning
problems, optimization algorithms are often run for a
limited number of iterations due to computational
constraints. In such settings, averaging schemes that
adapt to the transient and stationary phases of the
optimization process can lead to improved performance.

Furthermore, the interaction between learning rate
schedules and averaging methods plays a critical role
in determining the effectiveness of stochastic
optimization algorithms. Understanding these interactions
remains an important direction for both theoretical and
empirical research.

\section{Guidelines for Practitioners}
\label{sec:guidelines}

Averaging techniques are widely used in stochastic optimization,
yet the choice of an appropriate averaging strategy often depends
on the characteristics of the learning problem and the optimization
procedure. Based on the theoretical and empirical results reviewed
in this survey, we summarize several practical guidelines for the
use of averaging methods in stochastic gradient optimization.

\paragraph{When to Use Polyak--Ruppert Averaging}

Polyak--Ruppert averaging is particularly well suited for problems
where the objective function satisfies standard assumptions such
as smoothness and convexity. In such settings, uniform averaging
of stochastic gradient iterates provides strong theoretical
guarantees, including asymptotic optimality and variance reduction.
For statistical estimation problems and classical stochastic
approximation tasks, Polyak--Ruppert averaging remains one of
the most reliable approaches.

\paragraph{Handling the Transient Phase}

In many practical optimization problems, the early stages of
training are dominated by large updates as the algorithm moves
toward the vicinity of the optimum. Averaging these early
iterates may introduce bias into the final estimate. In such
cases, tail averaging or window averaging methods are often
preferable. These methods average only the most recent iterates
and therefore reduce the influence of the initial transient
phase of the optimization process.

\paragraph{Deep Learning Applications}

In deep learning applications, parameter averaging techniques
have proven particularly effective. Exponential moving averages
are widely used to stabilize training, while stochastic weight
averaging can improve generalization by identifying wider minima
of the loss landscape. These techniques are often easy to
implement and incur little additional computational cost,
making them attractive for large-scale neural network training.

\paragraph{Choice of Averaging Window}

The choice of averaging window or weighting scheme can have a
significant impact on performance. Larger averaging windows
generally provide stronger variance reduction, while smaller
windows may adapt more quickly to changes in the optimization
trajectory. In practice, the optimal window size often depends
on factors such as learning rate schedules, dataset size,
and the level of stochastic noise.

\paragraph{Computational Considerations}

Most averaging techniques are computationally inexpensive and
can be implemented with minimal additional memory. Uniform
averaging requires storing only a running sum of iterates,
while exponential moving averages can be updated recursively
with constant memory. These properties make averaging methods
particularly attractive in large-scale machine learning
applications where computational efficiency is essential.

Overall, averaging provides a simple yet powerful mechanism
for improving the stability and performance of stochastic
optimization algorithms. Practitioners are encouraged to
experiment with different averaging strategies to determine
which approach best suits their specific application.

\section{Open Problems and Future Directions}
\label{sec:open}

Although averaging techniques have been widely studied in stochastic
optimization, several important questions remain open. As stochastic
gradient methods continue to play a central role in machine learning,
understanding the theoretical and practical implications of averaging
remains an active area of research.

\paragraph{Finite-Sample Optimal Averaging}

Classical results on Polyak--Ruppert averaging focus primarily on
asymptotic optimality. In many machine learning applications,
however, optimization algorithms are run for a limited number of
iterations, and asymptotic guarantees may not accurately describe
practical behavior. Developing averaging schemes that are
\emph{optimal in finite-sample regimes} remains an important
open problem \cite{Schmidt2017}.

In particular, determining optimal weighting strategies for
finite sequences of stochastic gradient iterates is still
poorly understood. Such results could provide improved
guidelines for designing averaging schemes that adapt to the
transient and stationary phases of stochastic optimization.

\paragraph{Adaptive Averaging Strategies}

Most existing averaging methods use fixed rules for combining
iterates, such as uniform averaging or exponentially decaying
weights. However, the behavior of stochastic optimization
algorithms often varies across different phases of training.

Some adaptive averaging strategies have been studied \cite{Duchi2011,Wilson2017}. But there is a need to  automatically detect
changes in the optimization dynamics may provide improved
performance. For example, identifying when the optimization
process enters a stationary regime could allow averaging
methods to discard early transient iterates and focus on
later iterates that better reflect the local structure of
the objective function.

Developing theoretically justified adaptive averaging
schemes remains an open research direction.

\paragraph{Averaging in Deep Learning Optimization}

While averaging techniques have demonstrated strong empirical
performance in deep learning \cite{Goodfellow2016,Keskar2017,Sutskever2013}, a complete theoretical
understanding of their success remains limited. In particular,
it is still not fully understood why parameter averaging often
leads to improved generalization performance in neural networks.

Recent studies suggest that averaging may help identify
solutions located in flatter regions of the loss landscape.
However, the precise relationship between averaging,
optimization dynamics, and generalization remains an
important topic for future investigation.

\paragraph{Distributed and Federated Optimization}

As machine learning systems increasingly rely on distributed
computing environments, averaging has become a key component
of many large-scale training algorithms. In distributed \cite{verbraeken2020survey} and
federated learning systems \cite{zhang2021survey}, model parameters are frequently
aggregated across multiple workers using averaging
operations.

Understanding the interaction between local stochastic
optimization and global averaging procedures remains a
challenging problem. Issues such as communication delays,
heterogeneous data distributions, and asynchronous updates
introduce new complexities that require further study.

\paragraph{Connections with Other Variance Reduction Methods}

Averaging techniques are closely related to other approaches
for reducing variance in stochastic optimization algorithms,
including variance-reduced gradient methods and adaptive
optimization techniques. Exploring the connections between
these approaches may lead to new algorithmic insights and
improved optimization strategies. Developing unified frameworks that integrate averaging
methods with other variance reduction techniques represents
a promising direction for future research.

\bigskip

In summary, averaging remains a powerful and versatile tool
in stochastic optimization. From its origins in stochastic
approximation theory to its widespread use in modern machine
learning systems, averaging has played a crucial role in
improving the efficiency and stability of stochastic
gradient methods. Continued research on averaging
techniques is likely to yield further advances in both
theoretical understanding and practical optimization
algorithms.

\section{Conclusion}

Averaging techniques have played a central role in the development
of stochastic optimization methods. Beginning with the classical
stochastic approximation framework, iterate averaging has been
shown to significantly improve the stability and statistical
efficiency of stochastic gradient algorithms. Polyak--Ruppert
averaging established the theoretical foundation for these
methods, demonstrating that simple averaging of stochastic
iterates can achieve optimal asymptotic variance.

In modern machine learning, averaging methods have gained renewed
importance due to their ability to improve generalization and
stabilize large-scale training procedures. Techniques such as
exponential moving averages and stochastic weight averaging are
now widely used in deep learning systems. At the same time,
recent theoretical developments have improved our understanding
of the finite-sample behavior of averaging schemes and their
interaction with stochastic optimization dynamics.

Despite these advances, several challenges remain open. In
particular, designing adaptive averaging strategies and
understanding the role of averaging in nonconvex optimization
and deep learning remain active areas of research. As machine
learning models continue to grow in scale and complexity,
averaging techniques are likely to remain an important tool
for improving the efficiency and robustness of stochastic
optimization algorithms.

\bibliographystyle{plain}
\bibliography{avgsurv}

\end{document}